\documentclass[10pt,reqno]{amsart}
\usepackage{amsmath,amssymb,amsthm}
\usepackage{graphicx}
\usepackage[all]{xypic}
\input xypic
\def\thmsection{section}
\def\thmchangesection{changesection}

\def\thmchangechapter{changechapter}
\def\thmchange{change}
\def\thmplain{plain}
\ifx\theoremnumberstyle\thmplain
  \theoremstyle{break-italic}
  \newtheorem{satz}{Satz}
\else
  \ifx\theoremnumberstyle\thmsection
    \theoremstyle{break-italic}
    \newtheorem{satz}{Satz}[section]
  \else
    \ifx\theoremnumberstyle\thmchange
      \swapnumbers
      \theoremstyle{break-italic}
      \newtheorem{satz}{Satz}
    \else
      \ifx\theoremnumberstyle\thmchangesection
         \swapnumbers
         \theoremstyle{break-italic}
         \newtheorem{satz}{Satz}[section]
      \else
        \ifx\theoremnumberstyle\thmchangechapter
           \swapnumbers
           \theoremstyle{break-italic}
           \newtheorem{satz}{Satz}[chapter]
        \else
          \ifx\theoremnumberstyle\thmsection
             \theoremstyle{break-italic}
             \newtheorem{satz}{Satz}[section]
          \else
            \theoremstyle{break-italic}
            \newtheorem{satz}{Satz}[section]
          \fi
        \fi
      \fi
    \fi
  \fi
\fi

\theoremstyle{break-italic}
\newtheorem{theorem}[satz]{Theorem}
\newtheorem{lemma}[satz]{Lemma}

\newtheorem{corollary}[satz]{Corollary}
\newtheorem{Proposition}[satz]{Proposition}

\newtheorem*{conjecture*}{Conjecture}

\theoremstyle{break-roman}
\newtheorem{definition}[satz]{Definition}

\newtheorem{example}[satz]{Example}

\newtheorem{remark}[satz]{Remark}

\newtheorem{conjecture}[satz]{Conjecture}

\theoremstyle{standard}

\newtheorem*{claim}{Claim}

\theoremstyle{varthm-roman}
\newtheorem*{varthm-roman}{}
\theoremstyle{varthm-italic}
\newtheorem*{varthm-italic}{}
\theoremstyle{varthm-roman-break}
\newtheorem*{varthm-roman-break}{}
\theoremstyle{varthm-italic-break}
\newtheorem*{varthm-italic-break}{}
\theoremstyle{varthm-roman-no-punctuation}
\newtheorem{varthm-roman-no-punctuation-numbered}[satz]{}
\theoremstyle{varthm-italic-no-punctuation}
\newtheorem{varthm-italic-no-punctuation-numbered}[satz]{}

\newenvironment{varthm-roman-numbered}[1]{
  \begin{varthm-roman-no-punctuation-numbered}
    \mbox{\rm\textbf{#1}}
  }{\end{varthm-roman-no-punctuation-numbered}}
\newenvironment{varthm-italic-numbered}[1]{
  \begin{varthm-italic-no-punctuation-numbered}
    \mbox{\rm\textbf{#1}}
  }{\end{varthm-italic-no-punctuation-numbered}}
\newenvironment{varthm-roman-break-numbered}[1]{
  \begin{varthm-roman-no-punctuation-numbered}
    \mbox{\rm\textbf{#1}\newline}
  }{\end{varthm-roman-no-punctuation-numbered}}
\newenvironment{varthm-italic-break-numbered}[1]{
  \begin{varthm-italic-no-punctuation-numbered}
    \mbox{\rm\textbf{#1}}\newline
  }{\end{varthm-italic-no-punctuation-numbered}}
\numberwithin{equation}{section}
\def\nga{\widetilde}
\def\ex{\begin{example}
  }
  \def\eex{\end{example}}
\def\thr{\begin{theorem}}
\def\ethr{\end{theorem}}
\def\pro{\begin{Proposition}}
\def\epro{\end{Proposition}}
\def\coro{\begin{corollary}}
\def\ecoro{\end{corollary}}
\def\df{\begin{definition}}
\def\edf{\end{definition}}
\def\lm{\begin{lemma}}
\def\elm{\end{lemma}}
\def\pf{\begin{proof}}
\def\epf{\end{proof}}
\def\problem{\begin{problem}}
\def\eproblem{\end{problem}}

\def\ite{\begin{itemize}}
\def\hite{\end{itemize}}
\def\rem{\begin{remark}}
\def\erem{\end{remark}}

\def\cla{\begin{claim}}
\def\ecla{\end{claim}}
\def\conj{\begin{conjecture}}
\def\econj{\end{conjecture}}

\newcommand{\tru}{\setminus}
\newcommand{\td}{\Longleftrightarrow}
\newcommand{\sr}{\longrightarrow}
\newcommand{\Sr}{\Longrightarrow}

\newcommand{\mtn}{\rightarrow}

\newcommand{\gt}{\overline}

\newcommand{\seq}[1]{\left<#1\right>}

\def\Ann{\mbox{Ann}}
\def\depth{\mbox{depth}}

\def\Sing{\mbox{Sing}}
\def\NNor{\mbox{NNor}}

\def\NRed{\mbox{NRed}}
\def\Red{\mbox{Red}}

\def\c{\Bbb C}

\def\z{\Bbb Z}

\def\ohoa{\mathcal{O}}

\begin{document}
\title[Equinormalizability and Topologically triviality ] {Equinormalizability and topologically triviality of  deformations  of isolated curve singularities over smooth base spaces}

\author{ L\^{e} C\^{o}ng-Tr\`{i}nh}

\address{Department of Mathematics, Quy Nhon University, Vietnam }

\email{lecongtrinh@qnu.edu.vn}


\dedicatory{Dedicated to Professor Gert-Martin Greuel  on the occasion of his 70th birthday}
\begin{abstract} We give a $\delta$-constant  criterion for equinormalizability of  deformations   of isolated (not necessarily reduced) curve singularities over  smooth base spaces of dimension $\geq 1$.  For one-parametric families of isolated curve singularities, we show that their topologically triviality is equivalent to the admission of   weak simultaneous resolutions. 

\end{abstract}
\maketitle

\section{Introduction}
The theory of equinormalizable deformations has been initiated by B. Teissier  (\cite{Tei1})   in  the late 1970's for deformations of reduced curve singularities over $(\c,0)$.  It is generalized to higher dimensional base spaces by  M. Raynaud  and Teissier himself (\cite{Tei2}; some insight into the background of Raynaud's argument might be gleaned from the introduction to \cite{GrS}).  Recently, it is developed by Chiang-Hsieh and Lipman (\cite{Ch-Li}, 2006)  for projective deformations of reduced complex spaces over normal base spaces, and it is studied by  Koll\'{a}r (\cite{Ko}, 2011)  for projective deformations of generically reduced algebraic schemes over semi-normal base spaces.   

Each reduced curve singularity is associated with a $\delta$ number (see Definition \ref{dn4.2}), which is a finite number  and it is a topological invariant of  reduced curve singularities. Teissier-Raynaud-Chiang-Hsieh-Lipman (\cite{Tei1}, \cite{Tei2}, \cite{Ch-Li})  showed that a deformation of a reduced curve singularity over a normal base space  is equinormalizable (see Definition \ref{dn4.1}) if and only if it is $\delta$-constant, that is  the $\delta$ number of all of its fibers  are the same. This is so-called the \textit{$\delta$-constant criterion} for equinormalizability of deformations of reduced curve singularities. 

For   isolated curve singularities   with  embedded components,  Br\"{u}cker and Greuel   (\cite{BG}, 1990)   gave a similar $\delta$-constant criterion (with a new definition of the $\delta$ number, see Definition \ref{dn4.2}) for  equinormalizability of deformations of isolated (not necessarily reduced) curve singularities over $(\c,0)$.  The author considered  in \cite{Le} (2012)  deformations of \textit{plane curve singularities} with embedded components over smooth base spaces of dimension $\geq 1$, and gave a similar $\delta$-constant criterion for equinormalizability of these deformations, using special techniques (e.g. a corollary of Hilbert-Burch theorem), which are effective only for plane curve singularities. 

The first purpose of this paper is to generalize the $\delta$-constant criterion given in \cite{BG}  and \cite{Le} to deformations of isolated (not necessarily reduced)  curve singularities over normal or smooth base spaces of dimension $\geq 1$.  In Proposition  \ref{pro4.1} we show that equinormalizability of deformations of isolated curve singularities over normal base spaces implies the constancy of the $\delta$ number of  fibers of these deformations. Moreover, in  Theorem \ref{thr4.1} we show that if the normalization of the total space of 
a deformation of an isolated curve singularity over $(\c^{k},0)$, $k\geq 1$, is Cohen-Macaulay then the converse holds.  The assumption on Cohen-Macaulayness of the normalization of the total space ensures for flatness of the composition map. Moreover, Cohen-Macaulayness of the normalization of the total space is always satified for deformations over $(\c,0)$, because in this case, the total space is a   normal surface singularity, which is Cohen-Macaulay. 

In all of known results for  the $\delta$-constant criterion for  equinormalizability of deformations of isolated curve singularities, the total spaces of these deformations are always assumed to be reduced and pure dimensional.  It is necessary to weaken the hypothesis on   reducedness or purity of the dimension of total spaces. In section 2 we study the relationship between  reducedness of the total space and that of the \textit{generic fibers} of a flat morphism, and show in Theorem \ref{thr2.1}  that if the generic fibers of a flat morphism over a reduced Cohen-Macaulay space are reduced then the total space is reduced. In particular, if  there exists a representative of a deformation of an isolated  singularity over a reduced Cohen-Macaulay base space such that \textit{the total space is generically reduced over the base space} then the total space is reduced (see Corollary \ref{coro2.3}). This gives a way to check reducedness of the total space of a deformation, and to weaken the hypothesis on reducedness of the total space of a deformation. 

For families of  isolated curve singularities, one of the most important things is the admission of weak simultaneous resolutions (\cite{Tei2}) of these families. Buchweitz and Greuel (\cite{B-G}, 1980) gave a list of criteria  for the  admission of weak simultaneous resolutions of one-parametric families of  reduced curve singularities, namely,  the constancy of the Milnor number,  the constancy  of the $\delta$ number as well as   the number of branches of all  fibers, and the topologically triviality of these families (see Theorem \ref{thr5.1}). In the last section, we use a very new  result of Bobadilla, Snoussi and Spivakovsky (2014) to show that these criteria are also true for one-parametric families of isolated (not necessarily reduced) curve singularities  (see Theorem \ref{thr5.2}).     

\vspace{0.5cm}

\hspace{-0.6cm} \textbf{Notation:}  Let $f : (X,x) \mtn (S,0)$ be a morphism of complex  germs.   Denote by $(X^{red},x)$ the reduction of $(X,x)$ and $i: (X^{red},x) \hookrightarrow (X,x)$ the inclusion.  Let $\nu^{red}: (\gt{X},\gt{x}) \mtn (X^{red},x) $ be the normalization of   $(X^{red},x)$,  where $\gt{x}:=(\nu^{red})^{-1}(x)$.  Then the composition $\nu: (\gt{X},\gt{x})\overset{\nu^{red}}{\mtn} (X^{red},x) \overset{i}{\hookrightarrow} (X,x)$  is called the \textit{normalization of $(X,x)$}.  Denote $ \bar{f}:=f\circ \nu : (\gt{X}, \gt{x}) \mtn (S,0).$ For each $s\in S$, we denote 
$$ X_s:=f^{-1}(s), \quad \gt{X}_s:=\bar{f}^{-1}(s). $$

 \section{Generic reducedness}
Let  $f: (X,x)\mtn (S,0)$ be a flat morphism of complex germs. In this section we study  the relationship  between  reducedness of the total space $(X,x)$ and that of the  generic fibers of $f$.  This gives a way to check  reducedness of the total space of a flat morphism.

 \df Let $f: X \mtn S $ be a morphism of complex spaces. Denote by $\Red(X)$ the set of all reduced points of $X$ and
$$ \Red(f) = \{x \in X|  f \mbox{ is flat at } x \mbox{ and } f^{-1}(f(x)) \mbox{ is  reduced at } x\} $$
the \emph{reduced locus} of $f$.  We say
\ite
\item[(1)] $X$ is  \emph{generically reduced} if $\Red(X)$ is open and dense in $X$;
\item[(2)] $X$ is \emph{generically  reduced over}  $S$ if there is an analytically open dense set $V$ in $S$ such that $f^{-1}(V)$ is contained in $\Red(X)$;
\item[(3)]  the \emph{generic fibers of $f$ are reduced} if there is an analytically open dense set $V$ in $S$ such that $X_s:=f^{-1}(s)$ is  reduced for all $s$ in $V$.
\hite
\edf

We show in the following that under  properness of the restriction of a flat morphism $f: (X,x) \mtn (S,0)$ to its non-reduced locus, the generically reducedness of $X$ over $S$ implies reducedness of the generic fibers of $f$.

\pro \label{pro2.4} Let $f: (X,x) \mtn (S,0)$ be  flat with $(S,0)$ reduced. Assume that  there is a representative $f: X \mtn S$ such that its   restriction    on the non-reduced locus $\NRed(f):= X \tru \Red(f)$  is proper and  $X$ is generically reduced over $S$. Then  the  generic fibers of $f$  are reduced.
\epro
\pf    $\NRed(f)$ is analytically closed in $X$ (cf. \cite[Corollary I.1.116]{GLS}).  Moreover, since $X$ is generically reduced over $S$, there exists an analytically  open  dense set $U$ in $S$ such that $f^{-1}(U) \subseteq \Red(X)$.  Then, by properness of   the restriction $\NRed(f) \mtn S$,    $f(\NRed(f))$ is analytically closed and nowhere dense in $S$ by  \cite[Theorem 2.1(3), p.56]{BF}.   This implies that  $V:=S\tru f(\NRed(f))$ is analytically open dense in  $S$, and for  all $s \in  V$, $X_s : = f^{-1}(s)$ is  reduced.  Therefore the generic fibers of $f$ are reduced.
\epf

\coro  \label{coro2.5} Let $f: (X,x) \mtn (S,0)$ be   flat with $(S,0)$     reduced. Assume that  $X_0\tru \{x\}$ is reduced and  there exists a representative $f: X \mtn S$ such that $X$ is generically  reduced over $S$.  Then  the  generic fibers of $f$  are reduced. \\
In particular, if  $X_0 \tru \{x\}$ and $(X,x)$ are  reduced, then the generic fibers of $f$ are reduced.
\ecoro
\pf  Since $f$ is flat,  we have
$$ \NRed(f) \cap X_0 = \NRed(X_0) \subseteq \{x\},  $$
where $\NRed(X_0)$ denotes the set of non-reduced points of $X_0$.   This implies that the restriction $f: \NRed(f) \mtn S$ is finite, hence proper. Then the first assertion follows from Proposition \ref{pro2.4}.  Moreover, if $(X,x)$ is reduced then  there exists a representative $X$ of $(X,x)$ which  is reduced. Then  $X$ is obviously  generically reduced over some representative $S$ of $(S,s)$. Hence we have the latter assertion.
\epf

\rem \label{rem2.1} \rm The assumption on  reducedness of $X_0 \tru \{x\}$ in Corollary  \ref{coro2.5} is necessary  for reducedness of generic fibers, even for the case $S=\c$.  In fact, let $(X_0,0)\subseteq (\c^{3},0)$ be defined by the ideal
$$I_0=\seq{x^{2},y}\cap \seq{y^{2},z}\cap \seq{z^{2},x}\subseteq \c\{x,y,z\}$$
 and $(X,0)\subseteq (\c^{4},0)$ defined by the ideal
$$I=\seq{x^{2}-t^{2},y}\cap \seq{y^{2}-t^{2},z} \cap \seq{z^{2},x}\subseteq \c\{x,y,z,t\}.$$
 Let $f: (X,0) \mtn (\c,0)$ be the restriction on $(X,0)$ of the projection on the fourth component  $\pi: (\c^{4},0) \mtn (\c,0), ~(x,y,z,t)\mapsto t$. Then $f$ is flat, $X\tru X_0$ is reduced, hence $X$ is generically  reduced over some representative $T$ of $(\c,0)$. However  the fiber $(X_t,0)$ is not reduced for any $t\not = 0$.  Note that in this case $ X_0\tru \{0\}$ is not reduced.
\erem

 As we have seen from Corollary \ref{coro2.5},  if   the total space  of a  flat morphism over a reduced base space  is reduced, then   the  generic fibers of that morphism are reduced.  In the following we shows that over  a  reduced Cohen-Macaulay base space, the converse is also true. This generalizes  \cite[Proposition 3.1.1 (3)]{BG} to deformations over higher dimensional base spaces.

\thr   \label{thr2.1}  Let $f : (X,x) \mtn (S,0)$ be flat with  $(S,0)$  reduced   Cohen-Macaulay of dimension $k\geq 1$. If  there exists a representative $f: X \mtn S$ whose  generic fibers  are reduced,   then $(X,x)$  is reduced.
\ethr
\pf
 We divide the proof of this part into two steps.\\
\textbf{Step 1:} $\bold{S=\c^k}.$  Then  $f=(f_1,\cdots,f_k): (X,x) \mtn (\c^k,0)$ is flat. \\
For $k=1$, assume that there exists a representative $f: X \mtn T$ such that $X_t:=f^{-1}(t)$ is reduced for every $t\not =0$. Then for any $y \in X\tru X_0$ we have  $(X_{f(y)},y)$ is reduced. It follows that $(X,y)$ is reduced (cf. \cite[Theorem I. 1.101]{GLS}). Thus $X\tru X_0$ is reduced.   To show that $(X,x)$ is reduced, let $g$ be a nilpotent element of $\ohoa_{X,x}$. Then we have
$$ supp(g) = V(\Ann(g)) \subseteq X_0 = V(f).$$
 It follows from  Hilbert-R\"{u}ckert's  Nullstellensatz (cf. \cite[Theorem  I.1.72]{GLS}) that $f^n \in \Ann(g)$ for some $n\in \z_{+}$.
Hence $f^ng = 0$ in $\ohoa_{X,x}$.   Since $f$ is flat, it is a  non-zerodivisor in  $\ohoa_{X,x}$. Then $f^n$ is also a non-zerodivizor  in $\ohoa_{X,x}$. It follows that    $g=0$. Thus $(X,x)$ is reduced, and the statement is true for $k=1$.\\
For $k\geq 2$,  suppose  there is a representative $f: X \mtn S$  and an analytically open dense set $V$ in $S$ such that $X_s$ is  reduced for all $s\in V$. Let  us denote  by $H$ the line
$$H:= \{(t_1,\cdots,t_k) \in \c^k| t_1 = \cdots = t_{k-1}= 0 \}.$$
Denote by $A$ the complement of $V$ in $S$. Then $A$ is analytically  closed and nowhere dense in $S$.  We can choose coordinates $t_1,\cdots, t_k$ and a  representative of $(\c^k,0)$  such that $A \cap H = \{0\}$.  \\
    Denote $f':=(f_1,\cdots,f_{k-1})$.
Since $f$ is flat, $f_1,\cdots, f_{k-1}$ is an $\ohoa_{X,x}$-regular sequence, hence $f': (X,x) \mtn (\c^{k-1},0)$ is flat  with the special fiber  $(X',x): = (f'^{-1}(0),x) = (f^{-1}(H),x)$.    Since $f$ is flat,  $f_k$ is a non-zerodivisor in  $\ohoa_{X,x}/f'\ohoa_{X,x} = \ohoa_{X',x}$, hence  the morphism $f_k: (X',x) \mtn (\c,0)$ is flat.  For any $t\in \c\tru \{0\}$ close to $0$,
we have $(0,\cdots,0,t) \not \in  A$, hence   $f_k^{-1}(t) = f^{-1}(0,\cdots,0,t) $ is reduced. It follows from the case $k=1$  that the total space $(X',x)$ of $f_k$ is reduced. Since $f': (X,x) \mtn (\c^{k-1},0)$ is flat whose special fiber is  reduced,  $(X,x)$ is reduced (cf. \cite[Theorem I.1.101]{GLS}),  and we have the proof for this step. \\
{\bf Step 2:} $\bold{(S,0)}$ {\bf is Cohen-Macaulay of dimension } $\bold{k\geq 1}.$ Since $(S,0)$ is Cohen-Macaulay, there exists  an $\ohoa_{S,0}$-regular sequence $g_1, \cdots, g_k$, where $g_i \in \ohoa_{S,0} $ for every  $i = 1,\cdots, k$. Then  the  morphism
$$g=(g_1,\cdots,g_k): (S,0) \sr (\c^k,0),  t \longmapsto \big(g_1(t),\cdots,g_k(t)\big)$$
is flat.
We have
$$\dim  (g^{-1}(0),0)  = \dim \ohoa_{S,0}/(g_1,\cdots,g_k)\ohoa_{S,0} = 0$$ (cf. \cite[Prop. I.1.85]{GLS}).
This implies that  $g$ is finite.  Let $g: S \mtn T$ be  a  representative  which is flat and finite,  where $T$ is an open neighborhood of $0\in \c^k$.  Then   the composition   $h=g\circ f: X \sr T$  (for some representative) is  flat. To apply Step 1 for  $h$, we need to show the existence of an  analytically open  dense  set $U$  in  $T$ such that all fibers over $U$ are reduced. In fact,   since $S$ is  reduced,   its singular  locus  $ \Sing(S)$ is closed and nowhere dense  in $S$ (cf. \cite[Corollary I.1.111]{GLS}). It follows that  $A \cup  \Sing(S), $ $A$ as in Step 1,  is closed and nowhere dense  in $S$. Then the set     $U:=T\tru g(A\cup \Sing(S))$ is open and  dense in $T$ by the finiteness of $g$.  Furthermore, for any $t\in  U$,  $g^{-1}(t) = \{t_1,\cdots,t_r\}$, $t_i  \in V \cap (S\tru \Sing(S))$.  It follows that $h^{-1}(t) = f^{-1}(t_1) \cup \cdots \cup f^{-1}(t_r)$ is reduced. \\
Now applying   Step 1 for the flat map $h: X \mtn T$, we have reducedness of $(X,x)$.  The proof is complete.
 \epf
The following result  is a direct consequence of  Corollary \ref{coro2.5} and Theorem \ref{thr2.1}.
\coro \label{coro2.3} Let $f: (X,x) \mtn (S,0)$ be flat with $(S,0)$ reduced Cohen-Macaulay of dimension $k\geq 1$.  Suppose $X_0\tru \{x\}$ is reduced and there exists a representative $f: X \mtn S$ such that $X$ is generically reduced over $S$. Then $(X,x)$ is reduced.
\ecoro

Since normal surface singularities are reduced and Cohen-Macaulay, we have
\coro \label{coro2.4}\rm Let $f: (X,x) \mtn (S,0)$ be  flat with  $(S,0)$ a normal surface singularity. If  there exists a representative $f: X \mtn S$ whose  generic fibers  are reduced,   then $(X,x)$  is reduced.
 \ecoro

\section{Equinormalizable deformations of isolated curve singularities over smooth base spaces }

In this section we focus on equinormalizability of deformations of isolated (not necessarily reduced) curve singularities over smooth base spaces of dimension $\geq 1$. Because of isolatedness of singularities in the special fibers of these deformations, by Corollary \ref{coro2.3},  instead of assuming reducedness of the total spaces, we need only assume the generically reducedness of the total spaces over the base spaces.  

First we recall a definition of  equinormalizable deformations which follows Chiang-Hsieh-Lipman (\cite{Ch-Li}) and  Koll\'{a}r (\cite{Ko}). 
\df  \label{dn4.1}\rm    Let $f: X\sr S$ be a morphism of complex spaces. A \emph{simultaneous normalization of $f$ } is a morphism
$n: \nga{X} \sr X$ such that
 \ite
\item[(1)]  $n$ is finite,
\item[(2)] $\tilde{f}:=f\circ n: \nga{X}\mtn S$ is \emph{normal}, i.e.,  for each $z\in \nga{X}$, $\tilde{f}$ is flat at $z$ and the fiber $\nga{X}_{\tilde{f}(z)}:=\tilde{f}^{-1}(\tilde{f}(z))$ is normal,
\item[(3)] the induced map $n_s: \nga{X}_s:=\tilde{f}^{-1}(s) \sr X_s$ is bimeromorphic for each $s\in f(X)$.
\hite
The morphism $f$ is called \emph{equinormalizable} if  the normalization $\nu: \gt{X}\mtn X$ is a simultaneous normalization of $f$. It  is called \emph{ equinormalizable at $x\in X$} if the restriction of $f$ to some neighborhood of $x$ is equinormalizable.\\
If $f: (X,x) \sr (S,s)$ is a morphism  of germs, then a \emph{simultaneous normalization of $f$} is
a morphism $n$ from a multi-germ $(\nga{X}, n^{-1}(x))$ to $(X,x)$ such that some representative of $n$ is a simultaneous normalization of a representative of $f$. The germ $f$ is \emph{equinormalizable} if some representative of $f$ is equinormalizable.
 \edf
The following lemma allows us to do base change, reducing deformations over  higher dimensional base spaces to those over smooth  1-dimensional base spaces with similar properties.
\lm \label{lm4.1}
Let $f: (X,x) \mtn (S,0)$ be a deformation of an isolated singularity $(X_0,x)$ with $(S,0)$ normal.  Suppose that there exists some representative $f: X \mtn S$ such that $X$ is generically reduced over $S$.  Then  there exists  an open and dense set $U$ in  $S$ such that $X_s:=f^{-1}(s)$ is reduced, $\gt{X}_s:=\bar{f}^{-1}(s)$ is normal for all $s\in U$. Moreover,  for each $s\in U$, the induced morphism on the fibers $\nu_s:\gt{X}_s \mtn X_s$ is the normalization of $X_s$.
\elm
Here, we recall that $\nu: (\gt{X},\gt{x}) \mtn (X,x)$ is the normalization of $(X,x)$ and $\bar{f}:=f\circ \nu:  (\gt{X},\gt{x}) \mtn (S,0)$.
\pf  Since $X_0\tru \{x\}$ is reduced, it follows from the proof of Corollary \ref{coro2.5} that  the set $f(\NRed(f))$ is closed and nowhere dense in $S$. Denote by $\NNor(f)$ (resp. $\NNor(\bar{f})$)  the \textit{non-normal locus of  $f$ (resp. $\bar{f}$}), the set of points  $z$ in  $X$ (resp. $\gt{X}$) at which  either $f$ (resp. $\bar{f}$) is not flat or  $X_{f(z)}$ (resp. $\gt{X}_{\bar{f}(z)}$) is not normal.  Since $f$ is flat and $S$ is normal,  we   have $\nu(\NNor(\bar{f}) \cap \gt{X}_0) \subseteq \NNor(f) \cap X_0 = \NNor(X_0)$.  Equivalently, $\NNor(\bar{f}) \cap \gt{X}_0 \subseteq \nu^{-1}(\NNor(X_0))$ which is finite since $\nu$ is finite and $X_0$ has an isolated singularity at $x$.  It follows that the restriction of $\bar{f}$ on $\NNor(\bar{f})$ is finite. Then $\bar{f}(\NNor(\bar{f}))$ is closed and nowhere dense in $S$ by \cite[Theorem 2.1(3), p.56]{BF}. The set  $U:=S\tru \big(f(\NRed(f)) \cup \bar{f}(\NNor(\bar{f}))\big)$ satisfies all the required properties. 
\epf

 For deformations of isolated curve singularities we have the following necessary condition for their equinormalizability, in terms of the constancy of the $\delta$-invariant of fibers.  For the reader's convenience we recall the definition of the $\delta$-invariant of isolated (not necessarily reduced) curve singularities, which is defined by Br\"{u}cker and Greuel in \cite{BG}.

\df \label{dn4.2}  \rm Let $X$ be a complex curve and $x\in X$ an isolated singular point.  Denote by  $X^{red}$   its reduction and let  $\nu^{red}: \gt{X} \mtn X^{red}$  be  the normalization of the reduced curve $X^{red}$. The number
$$\delta(X^{red},x):=\dim_\c (\nu^{red}_*\ohoa_{\gt{X}})_x/\ohoa_{X^{red},x} $$
is called the \emph{delta-invariant of $X^{red}$ at $x$},
$$\epsilon(X,x):=\dim_\c H_{\{x\}}^0(\ohoa_X) $$
is called  the \emph{epsilon-invariant of $X$ at $x$},  where $H_{\{x\}}^0(\ohoa_X)$ denotes  local cohomology, and
$$\delta(X,x):=\delta(X^{red},x) - \epsilon(X,x) $$
is called the \emph{delta-invariant of $X$ at $x$}.\\
If $X$ has only finitely many singular points then the number
$$\delta(X):=\sum_{x\in \Sing(X)} \delta(X,x) $$
is called the \emph{delta-invariant } of $X$.
\edf
It is easy to see that $\delta(X^{red},x)\geq 0$, and $\delta(X^{red},x) = 0$ if  and only if $x$ is an isolated point of $X$ or the germ $(X^{red},x)$ is smooth.  Hence, if $x\in X$ is an isolated point of $X$ then $\delta(X,x) = -\dim_\c \ohoa_{X,x} = -  \epsilon(X,x)$. 
In particular,  $\delta(X,x) = -1$ for $x$ an isolated and reduced (hence normal) point of $X$.

\pro  \label{pro4.1}
Let $f: (X,x) \mtn (S,0)$  be a deformation of an isolated curve singularity $(X_0,x)$ with $(X,x)$  pure dimensional, $(S,0)$ normal. 
Suppose that there exists some representative $f: X \mtn S$ such that $X$ is generically reduced over $S$.   If $f$ is equinormalizable, then it is \textit{$\delta$-constant}, that is, $\delta(X_s) = \delta(X_0)$ for every $s\in S$ close to $0$.
\epro

\pf (Compare to the proof of \cite[Theorem 4.1 (2)]{Le})\\ 
It follows from Lemma \ref{lm4.1} that there exists an open and dense set $U$ in $S$ such that $X_s$ is reduced and $\gt{X}_s$ is normal for all $s\in U$. \\
We first show  that $f$ is $\delta$-constant on $U$,  i.e.  $\delta(X_s) = \delta(X_0)$ for any $s \in U$.   In fact, for any $s\in U$, $s\not =0$,  there exist an irreducible reduced curve singularity $C \subseteq S$ passing through $0$ and $s$. Let $\alpha: T \sr C \subseteq S$ be the normalization of this curve singularity such that $\alpha(T\tru \{0\}) \subseteq U$, where $T\subseteq \c$ is a small disc with center
 at $0$. Denote
$$X_T:=X\times_S T, ~ \gt{X}_T:= \gt{X}\times_S T.$$
 Then  we have the following Cartesian diagram:
 $$\xymatrix@C=12pt@R=10pt@M=8pt{
&&\ar @{} [dr] |{\Box} \gt{X}_T \ar[r] \ar[d]_{\nu_T} \ar@/_2pc/[dd]_{\bar{f}_T} & \gt{X} \ar[d]^{\nu} \ar@/^2pc/[dd]^{\bar{f}}\\
&&\ar @{} [dr] |{\Box} X_T\ar[d]_{f_T} \ar[r]& X \ar[d]^f\\
&&T \ar[r] & S}$$
For any $t\in T, s = \alpha(t) \in S$, we have
\begin{equation}\label{equ4.1}
\ohoa_{(X_T)_t}:= \ohoa_{f_T^{-1}(t)} \cong \ohoa_{X_s}, ~\ohoa_{(\gt{X}_T)_t}:= \ohoa_{\bar{f}_T^{-1}(t)} \cong \ohoa_{\gt{X}_s}.
\end{equation}
Since $f$ is flat by hypothesis and  $\bar{f}$ is flat by equinormalizability, it follows from the  preservation of flatness under base change
 (cf. \cite[Prop. I. 1.87]{GLS}) that the induced morphisms $f_T$ and $ \bar{f}_T$  are flat over $T$.  Hence, it follows from equinormalizability of $f$ and (\ref{equ4.1}) that $f_T: X_T \mtn T$ is equinormalizable. \\
 For any $t\in T\tru \{0\}$, $s=\alpha(t)\in U$, hence $(X_T)_t \cong  X_s$ is reduced by the existence of $U$. It follows from Theorem \ref{thr2.1} that $X_T$ is reduced.  On the other hand,  since $X$ and $S$ are pure dimensional, all  fibers of $f$, hence of   $f_T$,  are  pure dimensional by the dimension formula (\cite[Lemma, p.156]{Fi}).  Then  $X_T$ is also  pure dimensional because $T$ is pure 1-dimensional.   Therefore   it  follows from \cite[Korollar 2.3.5]{BG} that $f_T: X_T \sr T$ is $\delta$-constant, hence $f: X \sr S$ is $\delta$-constant on $U$. \\
Let us now take  $s_0 \in S\tru U$. Since $U$ is dense in $S$, $s_0 \in S$,
there exists always a  point $s_1 \in U$ which is close to $s_0$.  It follows from the
semi-continuity of the $\delta$-function (cf. \cite[Lemma 4.2]{Le}) that
$$ \delta(X_0) \geq \delta(X_{s_0}) \geq \delta(X_{s_1}).$$
Moreover,  $\delta(X_0) = \delta(X_{s_1})$ as shown above.  It implies that $\delta(X_{s_0})=\delta(X_0)$.   Hence $f : X \sr S$ is $\delta$-constant.
\epf

\rem \label{rem4.1} \rm  The complex spaces $X_T$ and $\gt{X}_T$ appearing in the proof of Proposition \ref{pro4.1} have the following properties:
\ite
\item[(1)] $X_T$ is reduced; $\gt{X}_T$ is reduced if $\bar{f}_T$ is flat;
\item[(2)] they  have the same normalization $\nga{X_T}$;
\item[(3)]  fibers of the compositions $\nga{X_T} \overset{\mu_T}{\mtn} \gt{X}_T \overset{\bar{f}_T}{\mtn} T$ and $\nga{X_T} \overset{\theta_T}{\mtn} X_T \overset{f_T}{\mtn} T$ coincide.  
\hite
In fact, as we have seen in the proof of Proposition \ref{pro4.1}, $X_T$ is reduced. Moreover, if $\bar{f}_T$ is flat, since its generic fibers are reduced (actually normal), $\gt{X}_T$ is reduced by Theorem \ref{thr2.1}. Therefore we have (1).  \\
Now we show (2).  Since finiteness and surjectivity are preserved under base change, $\nu_T$ is finite and surjective.  Let us denote by $\mu_T:\nga{X_T}\mtn \gt{X}_T$ the normalization of $\gt{X}_T$.  Then the composition $\theta_T:=\mu_T \circ \nu_T$ is finite and surjective. \\
Denote $A:=\NNor(f_T)$. Since $X_T$ is reduced, $A$ is nowhere dense in  $X_T$. Moreover, since $\nu_T$ is finite and surjective, it follows from Ritt's lemma (cf. \cite[Chapter 5, \S 3, p.102]{GR}) that the preimage $A':=\nu_T^{-1}(A)$ is  nowhere dense in $\gt{X}_T$.   Furthermore,  for any $z\not \in A'$, $y=\nu_T(z) \not \in A$, hence the fiber $(X_T)_t$ resp. $ X_s$ is normal at $y$ resp. $\alpha_T(y)$, where $t=f_T(y), s=\alpha(t)$. Thus $(X,\alpha_T(y))\cong (\gt{X},\bar{\alpha}_T(z))$. It follows that $(X_T,y) \cong (\gt{X}_T,z)$. Therefore $\gt{X}_T\tru A' \cong X_T\tru A$. Then $(\mu_T\circ \nu_T)^{-1}(A)$ is nowhere dense in $\nga{X_T}$ and we have the isomorphism 
$$ \nga{X_T}\tru  (\mu_T\circ \nu_T)^{-1}(A) = \nga{X_T}\tru  \mu_T^{-1}(A') \cong \gt{X}_T\tru A' \cong X_T\tru A. $$
Therefore   $\theta_T$ is  bimeromorphic, whence it is the normalization of $X_T$. (3) is obvious.
\erem 
The following theorem is the main result of this section, which asserts that  under certain conditions, the $\delta$-criterion is sufficient for equinormalizability of  deformations of  isolated curve singularities  over smooth base spaces of dimension $\geq 1$. This gives a generalization of \cite[Korollar 2.3.5]{BG}.

\thr \label{thr4.1} Let $f: (X,x) \mtn (\c^{k},0)$, $k\geq 1$, be a deformation of an isolated curve singularity $(X_0,x)$ with $(X,x)$ pure dimensional. Suppose that there exists a representative $f: X \mtn S$ such that $X$ is generically reduced over $S$.
If the normalization   $\gt{X}$ of $X$  is Cohen-Macaulay \footnotemark \footnotetext{This holds always  for $k=1$, since normal surfaces are Cohen-Macaulay.}  and $f$ is $\delta$-constant, then   $f$ is equinormalizable.
\ethr
\pf   First we show that Cohen-Macaulayness of $\gt{X}$ implies flatness of  the composition $\bar{f}$. Since $\gt{X}$ is Cohen-Macaulay and $S$ is smooth, it is sufficient to check that the dimension formula holds for $\bar{f}$ (cf. \cite[Proposition, p.158]{Fi}). 
But it is always the case, since for any $z\in \nu^{-1}(x)$, we have 
\begin{align*}
\dim (\gt{X},z) &= \dim (X,x) = \dim (X_0,x) + k \quad \quad \mbox{by flatness of $f$}\\
&= \dim (\gt{X}_0,z) + k.
\end{align*}
The latter equality follows from finiteness and surjectivity of  $\nu_0:(\gt{X}_0,z) \mtn (X_0,x)$.

  Let $U\subseteq S$ be the open dense set with properties described as in Lemma \ref{lm4.1}. For any $s\in U$,  let $C\subseteq S$ be an irreducible reduced curve singularity passing through $s$ and $0$ such that $C\cap (S\tru U) = \{0\}$.  Let $\alpha: T \sr C \subseteq S$ be the normalization of this curve singularity such that $\alpha(T\tru \{0\}) \subseteq U$, where $T\subseteq \c$ is a small disc with center at $0$.  Denote $X_T$ and $\gt{X}_T$ as in the proof of Proposition \ref{pro4.1}. Then, since $\bar{f}$ is flat, it follows from Remark \ref{rem4.1} that $X_T$ and $\gt{X}_T$ are reduced and they have the same normalization $\nga{X}_T$.  Consider  the following Cartesian diagram:
$$\xymatrix@C=12pt@R=10pt@M=8pt{
&&\ar @{} [dr]  \nga{X_T} \ar[d]^{\mu_T} \ar@/_1pc/[dd]_{\theta_T}\ar@/_3pc/[ddd]_{\nga{f}_T} & \\
&&\ar @{} [dr] |{\Box} \gt{X}_T\ar[r]^{\bar{\alpha}_T} \ar[d]^{\nu_T} \ar@/_1pc/[dd]_{\bar{f}_T} & \gt{X} \ar[d]_\nu \ar@/^1pc/[dd]^{\bar{f}}\\
&&\ar @{} [dr] |{\Box} X_T\ar[d]^{f_T} \ar[r]^{\alpha_T}& X \ar[d]_f\\
&&T\ar[r]_\alpha & S}$$
 Since fibers of $f$ and $f_T$ are isomorphic,  $f_T$ is $\delta$-constant   and $X_T$ is pure dimensional. Then it follows from \cite[Korollar 2.3.5]{BG} that $f_T$ is equinormalizable.  Therefore, by definition,  for each $t\in T$,  $(\nga{X})_t:=(\nga{f}_T)^{-1}(t)$ is  normal, and it is the normalization of $(X_T)_t$. 

Let us consider the flat map $\bar{f}_T : \gt{X}_T \mtn T$  and consider the normalization  $\mu_T: \nga{X_T} \mtn \gt{X}_T$ of $\gt{X}_T$. It follows from \cite[Proposition 1.2.2]{BG}  that the composition  $\bar{f}_T\circ \mu_T : \nga{X_T}\mtn T$ is flat.  Moreover, by the same argument as given in Remark \ref{rem4.1}, we can show that $(X_T)_t$ and $(\gt{X}_T)_t$  have the same normalization for each $t\in T$. Hence the restriction on the fibers $(\nga{X})_t \mtn (\gt{X}_T)_t$ is the normalization. Thus  by definition, $\bar{f}_T$ is equinormalizable. 
Then  $\bar{f}_T$ is $\delta$-constant  by Proposition \ref{pro4.1} (or by \cite[Korollar 2.3.5]{BG}).  This implies that for any $t\in T\tru \{0\}$, we have
$$ \delta(\gt{X}_0) = \delta((\gt{X}_T)_0) = \delta ((\gt{X}_T)_t) = 0 \mbox{ (since } (\gt{X}_T)_t \mbox{ is normal).}$$

Now we show that $\gt{X}_0$ is reduced.  First we show that   $\nu(\NNor(\gt{X}_0)) \subseteq \NNor(X_0)$. In fact, if $y \not \in \NNor(X_0)$ then $X_0$ is normal at $y$. Since $f$ is flat and $S$ is normal at $0$, $X$ is normal at $y$ (cf. \cite[Theorem I.1.101]{GLS}). Therefore we have the isomorphism $ (\gt{X},z) \overset{\cong}{\sr} (X,y)$ for every $z\in \nu^{-1}(y)$.  It induces an isomorphism on the fibers $(\gt{X}_0,z) \overset{\cong}{\sr} (X_0,y)$, hence $\gt{X}_0$ is normal at every point $z\in \nu^{-1}(y)$.  It follows that $y\not \in \nu(\NNor(\gt{X}_0))$.\\
Then, for any  $z\in \NNor(\gt{X}_0)$,  since $\NNor(X_0)$ is nowhere dense in $X_0$, by Ritt's lemma (cf. \cite[Chapter 5, \S 3, 2, p.103]{GR}) and by the dimension formula (when $f$ is flat)   we have
\begin{align*}
&\dim (\nu(\NNor(\gt{X}_0)), \nu(z)) \leq \dim (\NNor(X_0),\nu(z))  < \dim (X_0,\nu(z))\\
&= \dim (X,\nu(z)) - \dim (S,0)  = \dim (\gt{X},z) - \dim (S,0) \leq \dim (\gt{X}_0,z).
\end{align*}
Furthermore, the restriction $\nu_0: \gt{X}_0\sr X_0$ is finite. Hence
$$\dim (\nu(\NNor(\gt{X}_0)), \nu(z)) =\dim (\NNor(\gt{X}_0),z) ~\mbox{(cf. \cite[Corollary, p.141]{Fi})}. $$
It follows that for any $z\in \NNor(\gt{X}_0)$ we have
 $\dim (\NNor(\gt{X}_0),z) < \dim (\gt{X}_0,z)$,  i.e., $\NNor(\gt{X}_0)$ is nowhere dense in $\gt{X}_0$ by Ritt's lemma. This implies that $\gt{X}_0$ is generically normal, whence generically reduced. \\
 Moreover, for each $z\in \nu^{-1}(x)$, since $\bar{f}$ is flat and $\dim (\gt{X},z) = \dim (X,x) = k+1$, we have 
$$ \depth (\ohoa_{\gt{X}_0,z}) = \depth(\ohoa_{\gt{X},z}) - k \geq (k+1) - k =1. $$
On the other hand, we have 
$$ \dim (\gt{X}_0,z) = \dim (\gt{X},z) - k = 1. $$
Hence $ \depth (\ohoa_{\gt{X}_0,z}) \geq 1 = \min \{1, \dim (\gt{X}_0,z)\}$, i.e. $\gt{X}_0$ satisfies $(S_1)$ at every point $z\in \nu^{-1}(x)$. This implies  that $\gt{X}_0$ is reduced at every point of $\nu^{-1}(x)$. Then $\gt{X}_0$ is normal, and it is the normalization of $X_0$. It follows  that $f$ is equinormalizable. The proof is complete.
\epf 
The following example illustrates our main theorem. 
\ex[{\cite{St}, cf. \cite[Example 4.2]{Le}}] \label{ex4.1}\rm Let us consider the curve singularity  $ (X_0,0)\subseteq (\c^4,0)$ defined by the ideal 
$$I_0:= \seq{x^2 - y^3,z,w} \cap \seq{x,y,w} \cap \seq{x,y,z,w^2} \subseteq \c\{x,y,z,w\}.$$
 The curve singularity $(X_0,0)$ is a union of a cusp
$C$ in the plane $z=w=0,$ a straight line $L = \{x = y = w = 0\}$ and
an embedded non-reduced point $O = (0,0,0,0)$. Now we consider the restriction $f: (X,0)\mtn (\c^2,0)$ of the projection $\pi:(\c^6,0)\mtn (\c^2,0), ~ (x,y,z,w,u,v)\mapsto (u,v),$ to the complex germ  $(X,0)$  defined by the ideal 
$$I=\seq{x^2-y^3+uy^2,z,w} \cap \seq{x,y,w-v}\subseteq \c\{x,y,z,w,u,v\}.$$
It is easy to check  that $f$ is flat, $f^{-1}(0,0) = (X_0,0)$, the total space $(X,0)$ is reduced and pure $3$-dimensional, with  two 3-dimensional irreducible components. 

We have $\delta((X_0)^{red}) = 2$, $\epsilon(X_0)=1$, hence $\delta(X_0)=1$. Moreover, for each $u,v\in \c\tru \{0\}$, we have 
$$\delta(X_{(u,v)}) = \delta((X_{(u,v)})^{red}) - \epsilon(X_{(u,v)})= 1-0=1; $$
$$\delta(X_{(u,0)})= 2-1=1;\quad  \delta(X_{(0,v)}) = 1-0 =1.$$
Hence $f$ is $\delta$-constant. 

Moreover,  the normalizations of the first  component $(X_1,0)$ and the second component  $(X_2,0)$ of $(X,0)$  are given respectively by
$$ \nu_1: (\c^3,0) \mtn (X_1,0), \quad (T_1,T_2,T_3) \mapsto (0,0,T_1,T_3,T_2,T_3) $$
and  
$$ \nu_2: (\c^3,0) \mtn (X_2,0), \quad (T_1,T_2,T_3) \mapsto (T_3^3+T_1T_3,T_3^2+T_1,0,0,T_1,T_2). $$
Hence the composition maps are given respectively by 
$$ \bar{f}_1:  (\c^3,0) \mtn (\c,0), \quad (T_1,T_2,T_3) \mapsto (T_2,T_3)$$
and 
$$ \bar{f}_2:  (\c^3,0) \mtn (\c,0), \quad (T_1,T_2,T_3) \mapsto (T_1,T_2).$$
On both  components, $\bar{f}$ is flat with normal fibers, hence $f$ is equinormalizable. Note that,  in this example, the normalization of $(X,0)$ is smooth. All the computation given above can be easily done by \textbf{SINGULAR} (\cite{DGPS}).
\eex 

\section{Topologically  triviality of one-parametric families of isolated curve singularities}
In this section we consider one-parametric families of isolated (not necessarily reduced) curve singularities and show that the topologically triviality of these families is equivalent to the admission of  weak simultaneous resolutions  (\cite{Tei2}). 

Let $f: (X,x) \mtn (\c,0)$ be a deformation of an isolated curve singularity $(X_0,x)$ with $(X,x)$  pure dimensional. Let $f: X \mtn T$ be a \textit{good representative} (in the sense of \cite[\S 2.1, p.248]{B-G}) such that $X$ is generically reduced over $T$. Then $X$ is reduced by Corollary \ref{coro2.3}. Let $\nu: \gt{X} \mtn X$ be the normalization of $X$. Denote $\bar{f}:=f\circ \nu: \gt{X} \mtn T$.
\df[{cf. \cite{BG}}] \rm 
\ite
\item[(1)] $f$ is said to be \textit{topologically trivial } if there is a homeomorphism $h: X \overset{\approx}{\mtn} X_0 \times T$ such that $f=\pi \circ h$, where $\pi: X_0 \times T \mtn T$ is the projection.
\item[(2)] Assume that $f$ admits a section $\sigma: T \mtn X$ such that $X_t\tru \sigma(t)$ is smooth for all $t\in T$. Then $f$  admits a \textit{weak simultaneous resolution} if $f$ is equinormalizable  and  
$$ \big(\nu^{-1}(\sigma(T))\big)^{red} \cong \big(\nu^{-1}(\sigma(0))\big)^{red} \times T \quad (\mbox{over $T$}).$$
\hite  
\edf 

\rem[{cf. \cite{Tei2}}] \label{rem5.1}\rm  $f$ admits a weak simultaneous resolution  if and only if $f$ is equinormalizable and the number of branches $r(X_t,\sigma(t)) $ of $(X_t,\sigma(t))$ is constant for all $t\in T$.  
\erem 
Buchweitz and Greuel (1980)  proved the following result for  families of reduced curve singularities.
\thr[{\cite[Theorem 5.2.2]{B-G}}]  \label{thr5.1} Let $f: X \mtn T$ be a good representative of a flat family of reduced curves with section $\sigma: T\mtn X$ such that $X_t\tru \sigma(t)$ is smooth for each $t\in T$. Then the following conditions are equivalent:
\ite 
\item[(1)] $f$ admits a weak simultaneous resolution;
\item[(2)] the delta number $\delta(X_t,\sigma(t))$ and the number of branches  $r(X_t,\sigma(t))$ are constant for $t\in T$;
\item[(3)] the Milnor number $\mu(X_t, \sigma(t))$ is constant for $t\in T$;
\item[(4)] $f$ is topologically trivial.
\hite 
\ethr 
We shall show that this result is also true for families of isolated (not necessarily reduced) curve singularities.  Due to Br\"{u}cker and Greuel (\cite{BG}), we give  a new definition for the \textit{Milnor number} of a curve singulariy $C$ at an isolated singular  point $c\in C$, namely, 
  $$ \mu(C,c):= 2 \delta(C,c) - r(C,c) +1. $$
The Milnor number of $C$ is defined to be 
$$ \mu(C):=\sum_{c\in \Sing(C)} \mu(C,c). $$
To state and prove a similar result to Theorem \ref{thr5.1} we need the following result of Bobadilla, Snoussi and Spivakovsky (2014).
\lm[{\cite[Theorem 4.4]{BSS}}]  \label{lm5.1} Let $f: (X,x) \mtn (\c,0)$ be a deformation of an isolated curve singularity $(X_0,x)$ with $(X,x)$  reduced. Assume that the singular locus $\Sing(X,x)$ of $(X,x)$ is smooth of dimension 1. If $f$ is topologically trivial, then for any $z \in \nu^{-1}(x)$,  $\bar{f}: (\gt{X},z) \mtn (\c,0)$ is topologically trivial,  and the  normalization $(\gt{X},\nu^{-1}(x))$ of $(X,x)$ is smooth.
\elm  
The following theorem is the main result of this section.
\thr  \label{thr5.2} Let $f: (X,x) \mtn (\c,0)$ be a deformation of an isolated curve singularity $(X_0,x)$ with $(X,x)$  pure dimensional.
 Let $f: X \mtn T$ be a good representative with section $\sigma: T\mtn X$ such that $X_t\tru \sigma(t)$ is smooth for each $t\in T$ and 
 $X$ is generically reduced over $T$.   Assume that  $\Sing(X,x)$  is  smooth of dimension 1.  Then the following conditions are equivalent:
\ite 
\item[(1)] $f$ admits a weak simultaneous resolution;
\item[(2)] the delta  number $\delta(X_t,\sigma(t))$ and the number of branches  $r(X_t,\sigma(t))$ are constant for $t\in T$;
\item[(3)] the Milnor number $\mu(X_t, \sigma(t))$ is constant for $t\in T$;
\item[(4)] $f$ is topologically trivial.
\hite 
\ethr 
\pf  The equivalence of (1) and (2) follows from Theorem \ref{thr4.1} (for $k=1$) and Remark \ref{rem5.1}. (2) $\td (3)$ because of the definition of the Milnor number.  The implication $(1) \Sr (4)$ is proved by the same way for families of reduced curve singularities as given in the proof of the implication $(4) \Sr (6)$  of   \cite[Theorem 5.2.2]{B-G}. Now we prove that $(4) \Sr (1)$. 

  For convenience, let us assume that  $\nu^{-1}(x) = \{z_1,\cdots,z_r\}$.  Note that $\gt{X}_0:=\bar{f}^{-1}(0)$ is reduced, $\gt{X}_t:=\bar{f}^{-1}(t)$ is smooth for every  $t\not =0$ by \cite[Lemma 2.1.1]{BG}.  Therefore for every $i=1,\cdots, r$,   $\bar{f}: (\gt{X},z_i) \mtn (\c,0)$ is  a family of  reduced curve singularities with smooth general fibers, and  there exist  sections $\bar{\sigma}_1, \cdots, \bar{\sigma}_r: T \mtn \gt{X}$  such that   $\bar{\sigma}_i(0)=z_i$, $\nu^{-1}(\sigma(t)) =\{\bar{\sigma}_1(t), \cdots, \bar{\sigma}_r(t)\}$,  and  $\gt{X}_t\tru \bar{\sigma}_i(t)$ is smooth for every  $t\in T$ and for every $i=1,\cdots, r$.  \\
Assume that $f$ is topologically  trivial. Then it follows from Lemma \ref{lm5.1} that the deformation $\bar{f}: (\gt{X}, z_i)  \mtn (\c,0)$ of $(\gt{X}_0,z_i)$  is also topologically trivial for every $i=1,\cdots, r$.  Hence it follows from Theorem \ref{thr5.1}, applying  for  the flat  family  of reduced curve singularities $\bar{f}: (\gt{X},z_i)  \mtn (\c,0)$ with section  $\bar{\sigma}_i  : (\c,0) \mtn (\gt{X},z_i)$,  that the delta  number $\delta(\gt{X}_t,\bar{\sigma}_i(t))$ and the number of branches  $r(\gt{X}_t,\bar{\sigma}_i(t))$ are constant for $t\in T$. Then for $t\not =0$ we have 
$$ \delta(\gt{X}_0) = \delta(\gt{X}_t) = 0. $$
Hence $\gt{X}_0$ is normal. It follows that $f$ is equinormalizable.  On the other hand, the equinormalizability of $f$ over the smooth base space $(\c,0)$ implies that  for  every  $t\in T$ and for each $i=1,\cdots, r$,   the  induced map of $\nu$  on the fibers   $\nu_t: (\gt{X}_t, \bar{\sigma}_i(t)) \mtn (X_t,\sigma(t))$  is  the normalization of the corresponding  irreducible component of $(X_t,\sigma(t))$.   
It follows that the number of irreducible components of $(X_t, \sigma(t))$  is equal to the cardinality of $\nu^{-1}(\sigma(t))$, which is equal to $r$ for  every $t\in T$.  Hence $r(X_t,\sigma(t))$ is constant for every $t\in T$.  It follows that $f$ admits a weak simultaneous resolution, and we have (1). 
\epf 
\ex \rm Let us consider again the curve singularity  $ (X_0,0)\subseteq (\c^4,0)$ considered in Example \ref{ex4.1}  which is defined by the ideal 
$$I_0:= \seq{x^2 - y^3,z,w} \cap \seq{x,y,w} \cap \seq{x,y,z,w^2} \subseteq \c\{x,y,z,w\}.$$
  Now we consider the restriction $f: (X,0)\mtn (\c,0)$ of the projection $\pi:(\c^5,0)\mtn (\c,0), ~ (x,y,z,w,t)\mapsto t,$ to the complex germ  $(X,0)$  defined by the ideal 
$$I=\seq{x^2-y^3+ty^2,z,w} \cap \seq{x,y,w-t}\subseteq \c\{x,y,z,w,t\}.$$
We can check the following (all of them can be checked easily  by SINGULAR):
\ite 
\item[(1)] $f$ is flat;
\item[(2)] $(X,0)$ is reduced and pure $2$-dimensional, with two 2-dimensional irreducible components;
\item[(3)] $f$ is $\delta$-constant with $\delta(X_t) = 1$ for all $t\in \c$ close to $0$;
\item[(4)] $r(X_t)=2$ for all $t\in \c$ close to $0$;
\item[(5)] $f$ is equinormalizable;
\item[(6)] the normalization of each component of $(X,0)$ is $(\c^2,0)$, which is smooth.
\hite 
By Theorem \ref{thr5.2}, $f$ is topologically trivial.
\eex

\textbf{Acknowledgements.}  The author would like to express his gratitude to Professor Gert-Martin Greuel for his  valuable discussions, careful proof-reading and a lot of precise comments. He would also  like to  thank the anonymous referees for their  careful proof-reading and suggestions.   This research is funded by Vietnam National Foundation for Science and Technology Development 
(NAFOSTED) under the grant number 101.99-2013.24. This work is finished during the author's postdoctoral fellowship at the Vietnam 
Institute for Advanced Study in Mathematics (VIASM). He thanks VIASM for financial support and hospitality.


\end{document}